\def\sphinxdocclass{report}
\documentclass[letterpaper,12pt,english]{article}
\ifdefined\pdfpxdimen
   \let\sphinxpxdimen\pdfpxdimen\else\newdimen\sphinxpxdimen
\fi \sphinxpxdimen=49336sp\relax

\usepackage[margin=1in,marginparwidth=0.5in]{geometry}
\usepackage[utf8]{inputenc}
\ifdefined\DeclareUnicodeCharacter
  \DeclareUnicodeCharacter{00A0}{\nobreakspace}
\fi
\usepackage{cmap}
\usepackage[T1]{fontenc}
\usepackage{amsmath,amssymb,amstext}
\usepackage{babel}
\usepackage{times}
\usepackage[Bjarne]{fncychap}
\usepackage{longtable}
\usepackage{sphinx}

\usepackage{multirow}
\usepackage{eqparbox}

\usepackage{hyperref}
\usepackage{hypcap}
\addto\captionsenglish{}
\addto\captionsenglish{}
\addto\captionsenglish{\renewcommand{\literalblockname}{Listing }}

\addto\extrasenglish{}

  \usepackage{amssymb}
  \usepackage{rotating}
  \newcommand{\chapter}[1]{}  
  \newcommand{\ignore}[1]{}  

\title{COCO: The Large Scale Black-Box Optimization Benchmarking (bbob-largescale) Test Suite}
\date{Mar 26, 2019}
\release{0.9}
\date{\vspace{-1ex}}\author{
      Ouassim Ait Elhara$^1$,
      Konstantinos Varelas$^1$,
      Duc Manh Nguyen$^2$,
      Tea Tu\v{s}ar$^3$,\\
      Dimo Brockhoff$^1$,
      Nikolaus Hansen$^1$, 
      Anne Auger$^1$ 
  \\
    $^1$RandOpt team, Inria research centre Saclay and CMAP, Ecole Polytechnique, France
  \\
    $^2$Hanoi National University of Education, Vietnam
  \\
    $^3$Jo\v{z}ef Stefan Institute, Ljubljana, Slovenia
    }
\newcommand{\sphinxlogo}{}
\renewcommand{\releasename}{Release}
\makeindex

\begin{document}

\maketitle
\phantomsection\label{\detokenize{index::doc}}

\chapter{CHAPTERTITLE}
\label{\detokenize{index:coco-the-large-scale-black-box-optimization-benchmarking-bbob-largescale-test-suite}}\label{\detokenize{index:chaptertitle}}
\begin{abstract}
The \sphinxcode{bbob-largescale} test suite, containing 24 single-objective
functions in continuous domain, extends the well-known
single-objective noiseless \sphinxcode{bbob} test suite \phantomsection\label{\detokenize{index:id1}}{\hyperref[\detokenize{index:han2009}]{\sphinxcrossref{{[}HAN2009{]}}}}, which has been used since 2009 in
the \href{http://numbbo.github.io/workshops}{BBOB workshop series}, to large dimension. The core idea is to make the rotational
transformations \(\textbf{R}, \textbf{Q}\) in search space that
appear in the \sphinxcode{bbob} test suite computationally cheaper while retaining some desired
properties. This documentation presents an approach that replaces a full rotational transformation with a combination of a block-diagonal matrix and two permutation matrices in order to construct test functions whose computational and memory costs scale linearly in the dimension of the problem.
\end{abstract}\sphinxtableofcontents 
\newpage

\section{Introduction}
\label{\detokenize{index:introduction}}
In the \sphinxcode{bbob-largescale} test suite, we consider single-objective, unconstrained minimization problems
of the form
\begin{equation*}
\begin{split}\min_{x \in \mathbb{R}^n} \ f(x),\end{split}
\end{equation*}
with problem dimensions \(n \in \{20, 40, 80, 160, 320, 640\}.\)

The objective is to find, as quickly as possible, one or several solutions \(x\) in the search
space \(\mathbb{R}^n\) with \sphinxstyleemphasis{small} value(s) of \(f(x)\in\mathbb{R}\). We
generally measure the \sphinxstyleemphasis{time} of an optimization run as the number of calls to (queries of) the objective function \(f\).

We remind in the next sections some notations and definitions.

\subsection{Terminology}
\label{\detokenize{index:terminology}}\begin{description}
\item[{\sphinxstyleemphasis{function}}] \leavevmode
We talk about an objective \sphinxstyleemphasis{function} \(f\) as a parametrized mapping
\(\mathbb{R}^n\to\mathbb{R}\) with scalable input space, that is,
\(n\) is not (yet) determined. Functions are parametrized such that
different \sphinxstyleemphasis{instances} of the ``same'' function are available, e.g. translated
or rotated versions.

\item[{\sphinxstyleemphasis{problem}}] \leavevmode
We talk about a \sphinxstyleemphasis{problem}, \href{http://numbbo.github.io/coco-doc/C/coco\_8h.html\#a408ba01b98c78bf5be3df36562d99478}{\sphinxcode{coco\_problem\_t}}, as a specific \sphinxstyleemphasis{function
instance} on which an optimization algorithm is run. Specifically, a problem
can be described as the triple \sphinxcode{(dimension, function, instance)}. A problem
can be evaluated and returns an \(f\)-value. In the context of performance
assessment, a target \(f\)- or indicator-value is attached to each problem.
That is, a target value is added to the above triple to define a single problem
in this case.

\item[{\sphinxstyleemphasis{runtime}}] \leavevmode
We define \sphinxstyleemphasis{runtime}, or \sphinxstyleemphasis{run-length} as the \sphinxstyleemphasis{number of evaluations}
conducted on a given problem, also referred to as number of \sphinxstyleemphasis{function} evaluations.
Our central performance measure is the runtime until a given target value
is hit.

\item[{\sphinxstyleemphasis{suite}}] \leavevmode
A test- or benchmark-suite is a collection of problems, typically between
twenty and a hundred.

\end{description}

\subsection{Functions, Instances and Problems}
\label{\detokenize{index:functions-instances-and-problems}}
Each function is \sphinxstyleemphasis{parametrized} by the (input) dimension, \(n\), its identifier \(i\), and the instance number, \(j\),
that is:
\begin{equation*}
\begin{split}f_i^j \equiv f(n, i, j): \mathbb{R}^n \to \mathbb{R} \quad x \mapsto f_i^j (x) = f(n, i, j)(x).\end{split}
\end{equation*}
Varying \(n\) or \(j\) leads to a variation of the same function \(i\) of a given suite.
By fixing \(n\) and \(j\) for function \(f_i\), we define an optimization \sphinxstylestrong{problem}
\((n, i, j)\equiv(f_i, n, j)\) that can be presented to the optimization algorithm.
Each problem receives again an index in the suite, mapping the triple \((n, i, j)\) to a single
number.

We can think of \(j\) as an index to a continuous parameter vector setting,
as it parametrizes, among others things, translations and rotations. In
practice, \(j\) is the discrete identifier for single instantiations of
these parameters.

\subsection{Runtime and Target Values}
\label{\detokenize{index:runtime-and-target-values}}
In order to measure the runtime of an algorithm on a problem, we
establish a hitting time condition.
For a given problem \((f_i, n, j)\), we prescribe a \sphinxstylestrong{target value} \(t\) as a specific \(f\)-value
of interest \phantomsection\label{\detokenize{index:id2}}{\hyperref[\detokenize{index:han2016perf}]{\sphinxcrossref{{[}HAN2016perf{]}}}}.
For a single run, when an algorithm reaches or surpasses the target value \(t\)
on problem \((f_i, n, j)\), we say that it has \sphinxstyleemphasis{solved the problem} \((f_i, n, j, t)\) --- it was successful. %
\begin{footnote}[1]\sphinxAtStartFootnote
Note the use of the term \sphinxstyleemphasis{problem} in two meanings: as the problem the
algorithm is benchmarked on, \((f_i, n, j)\), and as the problem, \((f_i, n, j, t)\), an algorithm can
solve by hitting the target \(t\) with the runtime, \(\mathrm{RT}(f_i, n, j, t)\), or may fail to solve.
Each problem \((f_i, n, j)\) gives raise to a collection of dependent problems \((f_i, n, j, t)\).
Viewed as random variables, the events \(\mathrm{RT}(f_i, n, j, t)\) given \((f_i, n, j)\) are not
independent events for different values of \(t\).
\end{footnote}

The \sphinxstylestrong{runtime} is, then, the evaluation count when the target value \(t\) was
reached or surpassed for the first time.
That is, the runtime is the number of \(f\)-evaluations needed to solve the problem
\((f_i, n, j, t)\). %
\begin{footnote}[2]\sphinxAtStartFootnote
Target values are directly linked to a problem, leaving the burden to
properly define the targets with the designer of the benchmark suite.
The alternative is to present final \(f\)-values as results,
leaving the (rather unsurmountable) burden to interpret these values to the
reader.
Fortunately, there is an automatized generic way to generate target values
from observed runtimes, the so-called run-length based target values
\phantomsection\label{\detokenize{index:id8}}{\hyperref[\detokenize{index:han2016perf}]{\sphinxcrossref{{[}HAN2016perf{]}}}}.
\end{footnote}
\sphinxstyleemphasis{Measured runtimes are the only way how we assess the performance of an
algorithm.}
Observed success rates are generally translated into runtimes on a subset of
problems.

If an algorithm does not hit the target in a single run, its runtime remains
undefined --- while, then, this runtime is bounded from below by the number of evaluations
in this unsuccessful run.
The number of available runtime values depends on the budget the
algorithm has explored (the larger the budget, the more likely the target-values are reached).
Therefore, larger budgets are preferable --- however they should not come at
the expense of abandoning reasonable termination conditions. Instead,
restarts should be done \phantomsection\label{\detokenize{index:id5}}{\hyperref[\detokenize{index:han2016ex}]{\sphinxcrossref{{[}HAN2016ex{]}}}}.

\section{Overview of the Proposed \sphinxstyleliteralintitle{bbob-largescale} Test Suite}
\label{\detokenize{index:overview-of-the-proposed-bbob-largescale-test-suite}}
The \sphinxcode{bbob-largescale} test suite provides 24 functions in six dimensions (20, 40, 80, 160, 320 and 640) within
the \href{https://github.com/numbbo/coco}{COCO} framework \phantomsection\label{\detokenize{index:id9}}{\hyperref[\detokenize{index:han2016co}]{\sphinxcrossref{{[}HAN2016co{]}}}}. It is derived from the existing single-objective, unconstrained \sphinxcode{bbob} test suite with
modifications that allow the user to benchmark algorithms on high dimensional problems efficiently.
We will explain in this section how the \sphinxcode{bbob-largescale} test suite is built.

\subsection{The single-objective \sphinxstyleliteralintitle{bbob} functions}
\label{\detokenize{index:the-single-objective-bbob-functions}}
The \sphinxcode{bbob} test suite relies on the use of a number of raw functions from
which 24 \sphinxcode{bbob} functions are generated. Initially, so-called \sphinxstyleemphasis{raw} functions
are designed. Then, a series of transformations on these raw functions, such as
linear transformations (e.g., translation, rotation, scaling) and/or non-linear
transformations (e.g., \(T_{\text{osz}}, T_{\text{asy}}\))
will be applied to obtain the actual \sphinxcode{bbob} test functions. For example, the test function
\(f_{13}(\mathbf{x})\) (\href{http://coco.lri.fr/downloads/download15.03/bbobdocfunctions.pdf\#page=65}{Sharp Ridge function}) with (vector) variable \(\mathbf{x}\)
is derived from a raw function defined as follows:
\begin{equation*}
\begin{split}f_{\text{raw}}^{\text{Sharp Ridge}}(\mathbf{z}) = z_1^2 + 100\sqrt{\sum_{i=2}^{n}z_i^2}.\end{split}
\end{equation*}
Then one applies a sequence of transformations:
a translation by using the vector \(\mathbf{x}^{\text{opt}}\);
then a rotational transformation \(\mathbf{R}\); then a scaling transformation
\(\mathbf{\Lambda}^{10}\); then another rotational transformation \(\mathbf{Q}\)
to get the relationship
\(\mathbf{z} = \mathbf{Q}\mathbf{\Lambda}^{10}\mathbf{R}(\mathbf{x} - \mathbf{x}^{\text{opt}})\); and finally
a translation in objective space by using \(\mathbf{f}_{\text{opt}}\) to obtain the final
function in the testbed:
\begin{equation*}
\begin{split}f_{13}(\mathbf{x}) = f_{\text{raw}}^{\text{Sharp Ridge}}(\mathbf{z}) + \mathbf{f}_{\text{opt}}.\end{split}
\end{equation*}
There are two main reasons behind the use of transformations here:
\begin{enumerate}
\item {} 
provide non-trivial problems that cannot be solved by simply exploiting some of their properties (separability, optimum at fixed position, ...) and

\item {} 
allow to generate different instances, ideally of similar difficulty, of the same problem by using different (pseudo-)random transformations.

\end{enumerate}

Rotational transformations are used to avoid separability and thus coordinate system dependence in the test functions.
The rotational transformations consist in applying
an orthogonal matrix to the search space: \(x \rightarrow z = \textbf{R}x\), where \(\textbf{R}\) is the
orthogonal matrix.
While the other transformations used in the \sphinxcode{bbob} test suite could be naturally extended to
the large scale setting due to their linear complexity, rotational transformations have quadratic time and
space complexities. Thus, we need to reduce the complexity of these transformations in order for them to be usable, in practice, in the large scale setting.

\subsection{Extension to large scale setting}
\label{\detokenize{index:extension-to-large-scale-setting}}
Our objective is to construct a large scale test suite where the cost of a function call is
acceptable in higher dimensions while preserving the main characteristics of the original functions in the \sphinxcode{bbob}
test suite.
To this end, we will replace the full orthogonal matrices of the rotational transformations,
which would be too expensive in our large scale setting, with orthogonal transformations
that have linear complexity in the problem dimension: \sphinxstyleemphasis{permuted orthogonal block-diagonal matrices} (\phantomsection\label{\detokenize{index:id10}}{\hyperref[\detokenize{index:ait2016}]{\sphinxcrossref{{[}AIT2016{]}}}}).

Specifically, the matrix of a rotational transformation \(\textbf{R}\)
will be represented as:
\begin{equation*}
\begin{split}\textbf{R} = P_{\text{left}}BP_{\text{right}}.\end{split}
\end{equation*}
Here, \(P_{\text{left}} \text{ and } P_{\text{right}}\) are two permutation matrices %
\begin{footnote}[3]\sphinxAtStartFootnote
A \sphinxstyleemphasis{permutation matrix} is a square binary matrix that has exactly one entry of
1 in each row and each column and 0s elsewhere.
\end{footnote} and \(B\) is a
block-diagonal matrix of the form:
\begin{equation*}
\begin{split}B = \left(\begin{matrix}
B_1 & 0 & \dots & 0 \\
0 & B_2 & \dots & 0 \\
0 & 0 & \ddots & 0 \\
0 & 0 & \dots & B_{n_b}
\end{matrix}
\right),\end{split}
\end{equation*}
where \(n_b\) is the number of blocks and \(B_i, 1 \leq i \leq n_b\)
are square matrices of sizes \(s_i \times s_i\) satisfying \(s_i \geq 1\)
and \(\sum_{i=1}^{n_b}s_i = n\). In this case, the matrices
\(B_i, 1 \leq i \leq n_b\) are all orthogonal. Thus, the matrix \(B\)
is also an orthogonal matrix.

This representation allows the rotational transformation \(\textbf{R}\) to satisfy three
desired properties:
\begin{enumerate}
\item {} 
Have (almost) linear cost (due to the block structure of \(B\)).

\item {} 
Introduce non-separability.

\item {} 
Preserve the eigenvalues and therefore the condition number of the original function when it is convex quadratic (since \(\textbf{R}\) is orthogonal).

\end{enumerate}

\subsection{Generating the orthogonal block matrix \protect\(B\protect\)}
\label{\detokenize{index:generating-the-orthogonal-block-matrix}}
The block-matrices \(B_i, i=1,2,...,n_b\) will be uniformly distributed in the set of
orthogonal matrices of the same size. To this end, we first generate square matrices with
sizes \(s_i\) (\sphinxtitleref{i=1,2,...,n\_b}) whose entries are i.i.d. standard normally distributed.
Then we apply the Gram-Schmidt process to orthogonalize these matrices.

The parameters of this procedure include:
\begin{itemize}
\item {} 
the dimension of the problem \(n\),

\item {} 
the block sizes \(s_1, \dots, s_{n_b}\), where \(n_b\) is the number of blocks. In this test suite, we set \(s_i = s := \min\{n, 40\} \forall i=1,2,...,n_b\) (except, maybe, for the last block which can be smaller) %
\begin{footnote}[4]\sphinxAtStartFootnote
This setting allows to have the problems in dimensions 20 and 40 overlap between the \sphinxcode{bbob} test suite and its large-scale extension since in these dimensions, the block sizes coincide with the problem dimensions.
\end{footnote} and thus \(n_b = \lceil n/s \rceil\).

\end{itemize}

\subsection{Generating the permutation matrices \protect\(P\protect\)}
\label{\detokenize{index:generating-the-permutation-matrices}}
In order to generate the permutation matrix \(P\), we start from the identity matrix and apply, successively, a set of so-called \sphinxstyleemphasis{truncated uniform swaps}.
Each row/column (up to a maximum number of swaps) is swapped with a row/column chosen uniformly from the set of rows/columns within a fixed range \(r_s\).
A random order of the rows/columns is generated to avoid biases towards the first rows/columns.

Let \(i\) be the index of the first
variable/row/column to be swapped and \(j\) be the index of the second swap variable. Then
\begin{equation*}
\begin{split}j \sim U(\{l_b(i), l_b(i) + 1, \dots, u_b(i)\} \backslash \{i\}),\end{split}
\end{equation*}
where \(U(S)\) is the uniform distribution over the set \(S\) and \(l_b(i) = \max(1,i-r_s)\)
and \(l_b(i) = \min(n,i+r_s)\) with \(r_s\) a parameter of the approach.
If \(r_s \leq (n-1)/2\), the average distance between
the first and the second swap variable ranges from \((\sqrt{2}-1)r_s + 1/2\) (in the case of an
asymmetric choice for \(j\), i.e. when \(i\) is chosen closer to \(1\) or \(n\) than \(r_s\)) to
\(r_s/2 + 1/2\) (in the case of a symmetric choice for \(j\)). It is maximal when the first swap variable is at least \(r_s\)
away from both extremes or is one of them.

\sphinxstylestrong{Algorithm 1} below describes the process of generating a permutation using a
series of truncated uniform swaps with the following parameters:
\begin{itemize}
\item {} 
\(n\), the number of variables,

\item {} 
\(n_s\), the number of swaps.

\item {} 
\(r_s\), the swap range.

\end{itemize}

Starting with the identity permutation \(p\) and another permuation \(\pi\), drawn uniform
at random, we apply the swaps defined above
by taking \(p_{\pi}(1), p_{\pi}(2), \dots, p_{\pi}(n_s)\), successively, as
first swap variable. The resulting vector \(p\) will be the desired permutation.

\sphinxstyleemphasis{Algorithm 1: Truncated Uniform Permutations}
\begin{itemize}
\item {} 
Inputs: problem dimension \(n\), number of swaps \(n_s\), swap range \(r_s.\)

\item {} 
Output: a vector \(\textbf{p} \in \mathbb{N}^n\), defining a permutation.

\end{itemize}
\begin{enumerate}
\item {} 
\(\textbf{p} \leftarrow (1, \dots, n)\)

\item {} 
Generate a permutation \(\pi\) uniformly at random

\item {} 
\(\textbf{for } 1 \leq k \leq n_s \textbf{ do}\)

\item {} \begin{itemize}
\item {} 
\(i \leftarrow \pi(k)\), i.e., \(\textbf{p}_{\pi(k)}\) is the first swap variable

\end{itemize}

\item {} \begin{itemize}
\item {} 
\(l_b \leftarrow \max(1, i-r_s)\)

\end{itemize}

\item {} \begin{itemize}
\item {} 
\(u_b \leftarrow \min(n, i+r_s)\)

\end{itemize}

\item {} \begin{itemize}
\item {} 
\(S \leftarrow \{l_b, l_b + 1, \dots, u_b\} \backslash \{i\}\)

\end{itemize}

\item {} \begin{itemize}
\item {} 
Sample \(j\) uniformly at random in \(S\)

\end{itemize}

\item {} \begin{itemize}
\item {} 
Swap \(\textbf{p}_i\) and \(\textbf{p}_j\)

\end{itemize}

\item {} 
\(\textbf{end for}\)

\item {} 
\(\textbf{return p}\)

\end{enumerate}

In this test suite, we set \(n_s = n \text{ and } r_s = \lfloor n/3 \rfloor\). Some numerical
results in \phantomsection\label{\detokenize{index:id15}}{\hyperref[\detokenize{index:ait2016}]{\sphinxcrossref{{[}AIT2016{]}}}} show that with such parameters, the proportion of variables that are
moved from their original position when applying Algorithm 1 is approximately 100\% for all
dimensions 20, 40, 80, 160, 320, and 640 of the \sphinxcode{bbob-largescale} test suite.

\subsection{Implementation}
\label{\detokenize{index:implementation}}
Now, we describe how these changes to the rotational transformations are implemented
with the realizations of \(P_{\text{left}}BP_{\text{right}}\).
This will be illustrated through an example
on the Ellipsoidal function (rotated) \(f_{10}(\mathbf{x})\) (see the table in the next section), which is defined by
\begin{equation*}
\begin{split}f_{10}(\mathbf{x}) = \gamma(n) \times\sum_{i=1}^{n}10^{6\frac{i - 1}{n - 1}} z_i^2  + \mathbf{f}_{\text{opt}}, \text{with } \mathbf{z} = T_{\text{osz}} (\mathbf{R} (\mathbf{x} - \mathbf{x}^{\text{opt}})), \mathbf{R} = P_{1}BP_{2},\end{split}
\end{equation*}
as follows:

(i) First, we obtain the three matrices needed for the transformation, \(B, P_1, P_2\),
as follows:
\begin{quote}

\begin{sphinxVerbatim}[commandchars=\\\{\}]
\PYG{n}{coco\PYGZus{}compute\PYGZus{}blockrotation}\PYG{p}{(}\PYG{n}{B}\PYG{p}{,} \PYG{n}{seed1}\PYG{p}{,} \PYG{n}{n}\PYG{p}{,} \PYG{n}{s}\PYG{p}{,} \PYG{n}{n\PYGZus{}b}\PYG{p}{)}\PYG{p}{;}
\PYG{n}{coco\PYGZus{}compute\PYGZus{}truncated\PYGZus{}uniform\PYGZus{}swap\PYGZus{}permutation}\PYG{p}{(}\PYG{n}{P1}\PYG{p}{,} \PYG{n}{seed2}\PYG{p}{,} \PYG{n}{n}\PYG{p}{,} \PYG{n}{n\PYGZus{}s}\PYG{p}{,} \PYG{n}{r\PYGZus{}s}\PYG{p}{)}\PYG{p}{;}
\PYG{n}{coco\PYGZus{}compute\PYGZus{}truncated\PYGZus{}uniform\PYGZus{}swap\PYGZus{}permutation}\PYG{p}{(}\PYG{n}{P2}\PYG{p}{,} \PYG{n}{seed3}\PYG{p}{,} \PYG{n}{n}\PYG{p}{,} \PYG{n}{n\PYGZus{}s}\PYG{p}{,} \PYG{n}{r\PYGZus{}s}\PYG{p}{)}\PYG{p}{;}
\end{sphinxVerbatim}
\end{quote}
\begin{enumerate}
\setcounter{enumi}{1}
\item {} 
Then, whereever in the \sphinxcode{bbob} test suite, we use the following

\end{enumerate}
\begin{quote}

\begin{sphinxVerbatim}[commandchars=\\\{\}]
\PYG{n}{problem} \PYG{o}{=} \PYG{n}{transform\PYGZus{}vars\PYGZus{}affine}\PYG{p}{(}\PYG{n}{problem}\PYG{p}{,} \PYG{n}{R}\PYG{p}{,} \PYG{n}{b}\PYG{p}{,} \PYG{n}{n}\PYG{p}{)}\PYG{p}{;}
\end{sphinxVerbatim}

to make a rotational transformation, then in the \sphinxcode{bbob-largescale} test suite, we replace it with the three transformations

\begin{sphinxVerbatim}[commandchars=\\\{\}]
\PYG{n}{problem} \PYG{o}{=} \PYG{n}{transform\PYGZus{}vars\PYGZus{}permutation}\PYG{p}{(}\PYG{n}{problem}\PYG{p}{,} \PYG{n}{P2}\PYG{p}{,} \PYG{n}{n}\PYG{p}{)}\PYG{p}{;}
\PYG{n}{problem} \PYG{o}{=} \PYG{n}{transform\PYGZus{}vars\PYGZus{}blockrotation}\PYG{p}{(}\PYG{n}{problem}\PYG{p}{,} \PYG{n}{B}\PYG{p}{,} \PYG{n}{n}\PYG{p}{,} \PYG{n}{s}\PYG{p}{,} \PYG{n}{n\PYGZus{}b}\PYG{p}{)}\PYG{p}{;}
\PYG{n}{problem} \PYG{o}{=} \PYG{n}{transform\PYGZus{}vars\PYGZus{}permutation}\PYG{p}{(}\PYG{n}{problem}\PYG{p}{,} \PYG{n}{P1}\PYG{p}{,} \PYG{n}{n}\PYG{p}{)}\PYG{p}{;}
\end{sphinxVerbatim}
\end{quote}

Here, \(n\) is again the problem dimension, \(s\) the size of the blocks in \(B\), \(n_b:\)
the number of blocks, \(n_s:\) the number of swaps, and \(r_s:\) the swap range as presented previously.

\sphinxstylestrong{Important remark:} Although the complexity of \sphinxcode{bbob} test suite is reduced considerably by the above replacement of
rotational transformations, we recommend running the experiment on the \sphinxcode{bbob-largescale} test suite in parallel.

\section{Functions in \sphinxstyleliteralintitle{bbob-largescale} test suite}
\label{\detokenize{index:functions-in-bbob-largescale-test-suite}}
The table below presents the definition of all 24 functions of the \sphinxcode{bbob-largescale} test suite in detail. Beside the important
modification on rotational transformations, we also make two changes to the raw functions in the \sphinxcode{bbob} test suite.
\begin{itemize}
\item {} 
All functions, except for the Schwefel, Schaffer, Weierstrass, Gallagher, and Katsuura functions, are normalized by the parameter \(\gamma(n) = \min(1, 40/n)\) to have uniform target values that are comparable, in difficulty, over a wide range of dimensions.

\item {} 
The Discus, Bent Cigar and Sharp Ridge functions are generalized such that they have a constant proportion of distinct axes that remain consistent with the \sphinxcode{bbob} test suite.

\end{itemize}

For a better understanding of the properties of these functions and for the definitions
of the used transformations and abbreviations, we refer the reader to the original
\sphinxcode{bbob} \href{http://coco.lri.fr/downloads/download15.03/bbobdocfunctions.pdf}{function documention} for details.
\begin{sidewaystable}
    \centering
    \caption{Function descriptions of the separable, moderate, and ill-conditioned function groups of the {\ttfamily bbob-largescale} test suite.}
    \scriptsize
\noindent\begin{tabulary}{\linewidth}{|p{0.18 \textwidth}|p{0.41 \textwidth}|p{0.41 \textwidth}|}
\hline
\sphinxstylethead{\relax \unskip}\relax &\sphinxstylethead{\relax 
Formulation
\unskip}\relax &\sphinxstylethead{\relax 
Transformations
\unskip}\relax \\
\hline\multicolumn{3}{|l|}{\relax 
\sphinxstylestrong{Group 1: Separable functions}
\unskip}\relax \\
\hline
Sphere Function
&
\(f_1(\mathbf{x}) = \gamma(n) \times\sum_{i=1}^{n} z_i^2 + \mathbf{f}_{\text{opt}}\)
&
\(\mathbf{z} = \mathbf{x} - \mathbf{x}^{\text{opt}}\)
\\
\hline
Ellipsoidal Function
&
\(f_2(\mathbf{x}) = \gamma(n) \times\sum_{i=1}^{n}10^{6\frac{i - 1}{n - 1}} z_i^2+ \mathbf{f}_{\text{opt}}\)
&
\(\mathbf{z} = T_{\text{osz}}\left(\mathbf{x} - \mathbf{x}^{\text{opt}}\right)\)
\\
\hline
Rastrigin Function
&
\(f_3(\mathbf{x}) = \gamma(n) \times\left(10n - 10\sum_{i=1}^{n}\cos\left(2\pi z_i \right) + ||z||^2\right) + \mathbf{f}_{\text{opt}}\)
&
\(\mathbf{z} = \mathbf{\Lambda}^{10} T_{\text{asy}}^{0.2} \left( T_{\text{osz}}\left(\mathbf{x} - \mathbf{x}^{\text{opt}}\right) \right)\)
\\
\hline
Bueche-Rastrigin Function
&
\(f_4(\mathbf{x}) = \gamma(n) \times\left(10n - 10\sum_{i=1}^{n}\cos\left(2\pi z_i \right) + ||z||^2\right)\)
\(+ 100f_{pen}(\mathbf{x}) + \mathbf{f}_{\text{opt}}\)
&
\(z_i = s_i T_{\text{osz}}\left(x_i - x_i^{\text{opt}}\right) \text{for } i = 1,\dots, n\hspace{6cm}\)
\(s_i = \begin{cases} 10 \times 10^{\frac{1}{2} \frac{i-1}{n-1}} & \text{if } z_i >0 \text{ and } i \text{ odd} \\ 10^{\frac{1}{2} \frac{i-1}{n-1}} & \text{otherwise} \end{cases}\)
\(\text{ for } i = 1,\dots, n\)
\\
\hline
Linear Slope
&
\(f_5(\mathbf{x}) = \gamma(n)\times \sum_{i=1}^{n}\left( 5 \vert s_i \vert - s_i z_i \right) + \mathbf{f}_{\text{opt}}\)
&
\(z_i = \begin{cases} x_i & \text{if } x_i^{\mathrm{opt}}x_i < 5^2 \\ x_i^{\mathrm{opt}} & \text{otherwise} \end{cases}\)
\(\text{ for } i=1, \dots, n,\hspace{3.5cm}\)
\(s_i = \text{sign} \left(x_i^{\text{opt}}\right) 10^{\frac{i-1}{n-1}} \text{ for } i=1, \dots, n,\hspace{4cm}\)
\(\mathbf{x}^{\text{opt}} = \mathbf{z}^{\text{opt}} = 5\times \mathbf{1}_{-}^+\)
\\
\hline\end{tabulary}

\noindent\begin{tabulary}{\linewidth}{|p{0.18 \textwidth}|p{0.41 \textwidth}|p{0.41 \textwidth}|}
\hline
\multicolumn{3}{|l|}{\relax 
\sphinxstylestrong{Group 2: Functions with low or moderate conditioning}
\unskip}\relax \\
\hline
Attractive Sector Function
&
\(f_6(\mathbf{x}) = T_{\text{osz}}\left(\gamma(n) \times \sum_{i=1}^{n}\left( s_i z_i\right)^2 \right)^{0.9} + \mathbf{f}_{\text{opt}}\)
&
\(\mathbf{z} = \mathbf{Q} \mathbf{\Lambda}^{10} \mathbf{R}(\mathbf{x} - \mathbf{x}^{\text{opt}})\)
\(\hspace{0.2cm} \text{ with } \mathbf{R} = P_{11}B_1P_{12}, \mathbf{Q} = P_{21}B_2P_{22},\hspace{1.5cm}\)
\(s_i = \begin{cases} 10^2 & \text{if } z_i \times x_i^{\mathrm{opt}} > 0\\ 1 & \text{otherwise}\end{cases}\)
\(\text{ for } i=1,\dots, n\)
\\
\hline
Step Ellipsoidal Function
&
\(f_7(\mathbf{x}) = \gamma(n) \times 0.1 \max\left(\vert \hat{z}_1\vert/10^4, \sum_{i=1}^{n}10^{2\frac{i - 1}{n - 1}}z_i^2\right) + f_{pen}(\mathbf{x}) + \mathbf{f}_{\text{opt}}\)
&
\(\mathbf{\hat{z}} = \mathbf{\Lambda}^{10} \mathbf{R}(\mathbf{x}-\mathbf{x}^{\text{opt}})  \text{ with }\mathbf{R} = P_{11}B_1P_{12},\hspace{4.5cm}\)
\(\tilde{z}_i= \begin{cases} \lfloor 0.5 + \hat{z}_i \rfloor & \text{if }  |\hat{z}_i| > 0.5 \\ \lfloor 0.5 + 10 \hat{z}_i \rfloor /10 & \text{otherwise} \end{cases}\)
\(\text{ for } i=1,\dots, n,\hspace{1.5cm}\)
\(\mathbf{z} = \mathbf{Q} \mathbf{\tilde{z}} \text{ with } \mathbf{Q} = P_{21}B_2P_{22}\)
\\
\hline
Rosenbrock Function, original
&
\(f_8(\mathbf{x}) = \gamma(n) \times\sum_{i=1}^{n} \left(100 \left(z_{i}^2 - z_{i+1}\right)^2 + \left(z_{i} - 1\right)^2\right) + \mathbf{f}_{\text{opt}}\)
&
\(\mathbf{z} = \max\left(1, \dfrac{\sqrt{s}}{8}\right)(\mathbf{x} - \mathbf{x}^{\text{opt}})+ \mathbf{1},\)
\(\mathbf{x}^{\text{opt}} \in [-3, 3]^n\)
\\
\hline
Rosenbrock Function, rotated
&
\(f_9(\mathbf{x}) = \gamma(n) \times\sum_{i=1}^{n} \left(100 \left(z_{i}^2 - z_{i+1}\right)^2 + \left(z_{i} - 1\right)^2\right) + \mathbf{f}_{\text{opt}}\)
&
\(\mathbf{z} = \max\left(1, \dfrac{\sqrt{s}}{8}\right)\mathbf{R} (\mathbf{x} - \mathbf{x}^{\text{opt}})+ \mathbf{1}\)
\(\text{ with }\mathbf{R} = P_{1}BP_{2},\)
\(\mathbf{x}^{\text{opt}} \in [-3, 3]^n\)
\\
\hline\end{tabulary}

\noindent\begin{tabulary}{\linewidth}{|p{0.18 \textwidth}|p{0.41 \textwidth}|p{0.41 \textwidth}|}
\hline
\multicolumn{3}{|l|}{\relax 
\sphinxstylestrong{Group 3: Functions with high conditioning and unimodal}
\unskip}\relax \\
\hline
Ellipsoidal Function
&
\(f_{10}(\mathbf{x}) = \gamma(n) \times\sum_{i=1}^{n}10^{6\frac{i - 1}{n - 1}} z_i^2  + \mathbf{f}_{\text{opt}}\)
&
\(\mathbf{z} = T_{\text{osz}} (\mathbf{R} (\mathbf{x} - \mathbf{x}^{\text{opt}})) \text{ with }\mathbf{R} = P_{1}BP_{2}\)
\\
\hline
Discus Function
&
\(f_{11}(\mathbf{x}) = \gamma(n) \times\left(10^6\sum_{i=1}^{\lceil n/40 \rceil}z_i^2 + \sum_{i=\lceil n/40 \rceil+1}^{n}z_i^2\right) + \mathbf{f}_{\text{opt}}\)
&
\(\mathbf{z} = T_{\text{osz}}(\mathbf{R}(\mathbf{x} - \mathbf{x}^{\text{opt}})) \text{ with }\mathbf{R} = P_{1}BP_{2}\)
\\
\hline
Bent Cigar Function
&
\(f_{12}(\mathbf{x}) = \gamma(n) \times\left(\sum_{i=1}^{\lceil n/40 \rceil}z_i^2 + 10^6\sum_{i=\lceil n/40 \rceil + 1}^{n}z_i^2 \right) + \mathbf{f}_{\text{opt}}\)
&
\(\mathbf{z} = \mathbf{R} T_{\text{asy}}^{0.5}(\mathbf{R}((\mathbf{x} - \mathbf{x}^{\text{opt}})) \text{ with }\mathbf{R} = P_{1}BP_{2}\)
\\
\hline
Sharp Ridge Function
&
\(f_{13}(\mathbf{x}) = \gamma(n) \times\left(\sum_{i=1}^{\lceil n/40 \rceil}z_i^2 + 100\sqrt{\sum_{i=\lceil n/40 \rceil + 1}^{n}z_i^2} \right) + \mathbf{f}_{\text{opt}}\)
&
\(\mathbf{z} = \mathbf{Q}\mathbf{\Lambda}^{10}\mathbf{R}(\mathbf{x} - \mathbf{x}^{\text{opt}})\)
\(\text{ with } \mathbf{R} = P_{11}B_1P_{12}, \mathbf{Q} = P_{21}B_2P_{22}\)
\\
\hline
Different Powers Function
&
\(f_{14}(\mathbf{x}) = \gamma(n) \times\sum_{i=1}^{n} \vert z_i\vert ^{\left(2 + 4 \times \frac{i-1}{n- 1}\right)} + \mathbf{f}_{\text{opt}}\)
&
\(\mathbf{z} = \mathbf{R}(\mathbf{x} - \mathbf{x}^{\text{opt}}) \text{ with }\mathbf{R} = P_{1}BP_{2}\)
\\
\hline\end{tabulary}

\end{sidewaystable}

\begin{sidewaystable}
    \centering
    \caption{Function descriptions of the multi-modal function group with adequate global structure of the {\ttfamily bbob-largescale} test suite.}
    \scriptsize
\noindent\begin{tabulary}{\linewidth}{|p{0.18 \textwidth}|p{0.41 \textwidth}|p{0.41 \textwidth}|}
\hline
\sphinxstylethead{\relax \unskip}\relax &\sphinxstylethead{\relax 
Formulation
\unskip}\relax &\sphinxstylethead{\relax 
Transformations
\unskip}\relax \\
\hline\multicolumn{3}{|l|}{\relax 
\sphinxstylestrong{Group 4: Multi-modal functions with adequate global structure}
\unskip}\relax \\
\hline
Rastrigin Function
&
\(f_{15}(\mathbf{x}) = \gamma(n) \times\left(10n - 10\sum_{i=1}^{n}\cos\left(2\pi z_i \right) + ||\mathbf{z}||^2\right) + \mathbf{f}_{\text{opt}}\)
&
\(\mathbf{z} = \mathbf{R} \mathbf{\Lambda}^{10} \mathbf{Q} T_{\text{asy}}^{0.2} \left(T_{\text{osz}} \left(\mathbf{R}\left(\mathbf{x} - \mathbf{x}^{\text{opt}} \right) \right) \right) \hspace{5cm}\)
\(\text{with } \mathbf{R} = P_{11}B_1P_{12}, \mathbf{Q} = P_{21}B_2P_{22}\)
\\
\hline
Weierstrass Function
&
\(f_{16}(\mathbf{x}) = 10\left( \dfrac{1}{n} \sum_{i=1}^{n} \sum_{k=0}^{11} \dfrac{1}{2^k} \cos \left( 2\pi 3^k \left( z_i + 1/2\right) \right) - f_0\right)^3\)
\(+\dfrac{10}{n}f_{pen}(\mathbf{x}) + \mathbf{f}_{\text{opt}}\)
&
\(\mathbf{z} = \mathbf{R}\mathbf{\Lambda}^{1/100}\mathbf{Q}T_{\text{osz}}(\mathbf{R}(\mathbf{x} - \mathbf{x}^{\text{opt}}))\hspace{6cm}\)
\(\text{with } \mathbf{R} = P_{11}B_1P_{12}, \mathbf{Q} = P_{21}B_2P_{22},\hspace{5.8cm}\)
\(f_0= \sum_{k=0}^{11} \dfrac{1}{2^k} \cos(\pi 3^k)\)
\\
\hline
Schaffers F7 Function
&
\(f_{17}(\mathbf{x}) = \left(\dfrac{1}{n-1} \sum_{i=1}^{n-1} \left(\sqrt{s_i} + \sqrt{s_i}\sin^2\left( 50 (s_i)^{1/5}\right)\right)\right)^2\)
\(+ 10 f_{pen}(\mathbf{x}) + \mathbf{f}_{\text{opt}}\)
&
\(\mathbf{z} = \mathbf{\Lambda}^{10} \mathbf{Q} T_{\text{asy}}^{0.5}(\mathbf{R}(\mathbf{x} - \mathbf{x}^{\text{opt}}))\)
\(\text{with } \mathbf{R} = P_{11}B_1P_{12}, \mathbf{Q} = P_{21}B_2P_{22},\hspace{1cm}\)
\(s_i= \sqrt{z_i^2 + z_{i+1}^2}, i=1,\dots, n-1\)
\\
\hline
Schaffers F7 Function,
moderately ill-conditioned
&
\(f_{18}(\mathbf{x}) = \left(\dfrac{1}{n-1} \sum_{i=1}^{n-1} \left(\sqrt{s_i} + \sqrt{s_i}\sin^2\left( 50 (s_i)^{1/5}\right)\right)\right)^2\)
\(+ 10 f_{pen}(\mathbf{x}) + \mathbf{f}_{\text{opt}}\)
&
\(\mathbf{z} = \mathbf{\Lambda}^{1000} \mathbf{Q} T_{\text{asy}}^{0.5}(\mathbf{R}(\mathbf{x} - \mathbf{x}^{\text{opt}}))\)
\(\text{ with } \mathbf{R} = P_{11}B_1P_{12}, \mathbf{Q} = P_{21}B_2P_{22},\hspace{0.5cm}\)
\(s_i= \sqrt{z_i^2 + z_{i+1}^2}, i=1,\dots, n-1\)
\\
\hline
Composite Griewank-Rosenbrock
Function F8F2
&
\(f_{19}(\mathbf{x}) = \dfrac{10}{n-1} \sum_{i=1}^{n-1} \left( \dfrac{s_i}{4000} - \cos\left(s_i \right)\right) + 10 + \mathbf{f}_{\text{opt}}\)
&
\(\mathbf{z} = \max\left(1, \dfrac{\sqrt{s}}{8}\right)\mathbf{R} \mathbf{x} + \dfrac{\mathbf{1}}{2}\)
\(\text{ with }\mathbf{R} = P_{1}BP_{2},\hspace{3.4cm}\)
\(s_i= 100(z_i^2 - z_{i+1})^2 + (z_i - 1)^2,\)
\(\text{ for } i=1,\dots, n-1,\)
\(\mathbf{z}^{\text{opt}} = \mathbf{1}\)
\\
\hline\end{tabulary}

\end{sidewaystable}

\begin{sidewaystable}
    \centering
    \caption{Function descriptions of the multi-modal function group with weak global structure of the {\ttfamily bbob-largescale} test suite.}
    \scriptsize
\noindent\begin{tabulary}{\linewidth}{|p{0.18 \textwidth}|p{0.41 \textwidth}|p{0.41 \textwidth}|}
\hline
\sphinxstylethead{\relax \unskip}\relax &\sphinxstylethead{\relax 
Formulation
\unskip}\relax &\sphinxstylethead{\relax 
Transformations
\unskip}\relax \\
\hline\multicolumn{3}{|l|}{\relax 
\sphinxstylestrong{Group 5: Multi-modal functions with weak global structure}
\unskip}\relax \\
\hline
Schwefel Function
&
\(f_{20}(\mathbf{x}) = -\dfrac{1}{100 n} \sum_{i=1}^{n} z_i\sin\left(\sqrt{\vert z_i\vert}\right) + 4.189828872724339\)
\(+ 100f_{pen}(\mathbf{z}/100)+\mathbf{f}_{\text{opt}}\)
&
\(\mathbf{\hat{x}} = 2 \times \mathbf{1}_{-}^{+} \otimes \mathbf{x},\)
\(\hat{z}_1 = \hat{x}_1, \hat{z}_{i+1}=\hat{x}_{i+1} + 0.25 \left(\hat{x}_{i} - 2\left|x_i^{\text{opt}}\right|\right),\)
\(\text{ for } i=1, \dots, n-1,\)
\(\mathbf{z} = 100 \left(\mathbf{\Lambda}^{10} \left(\mathbf{\hat{z}} - 2\left|\mathbf{x}^{\text{opt}}\right|\right) + 2\left|\mathbf{x}^{\text{opt}}\right|\right),\)
\(\mathbf{x}^{\text{opt}} = 4.2096874633/2 \mathbf{1}_{-}^{+}\)
\\
\hline
Gallagher's Gaussian
101-me Peaks Function
&
\(f_{21}(\mathbf{x}) = T_{\text{osz}}\left(10 - \max_{i=1}^{101} w_i \exp\left(- \dfrac{1}{2n} (\mathbf{z} - \mathbf{y}_i)^T\mathbf{B}^T\mathbf{C_i}\mathbf{B} (\mathbf{z} - \mathbf{y}_i) \right) \right)^2\)
\(+ f_{pen}(\mathbf{x}) + \mathbf{f}_{\text{opt}}\)
&
\(w_i = \begin{cases} 1.1 + 8 \times \dfrac{i-2}{99} & \text{for } 2 \leq i \leq 101\\ 10 & \text{for } i = 1 \end{cases}\)
\(\mathbf{B} \text{ is a block-diagonal matrix without}\)
\(\text{permuations of the variables.}\)
\(\mathbf{C_i} = \Lambda^{\alpha_i}/\alpha_i^{1/4} \text{where } \Lambda^{\alpha_i} \text{ is defined as usual,}\)
\(\text{but with randomly permuted diagonal elements.}\)
\(\text{For } i=2,\dots, 101, \alpha_i \text{ is drawn uniformly}\)
\(\text{from the set } \left\{1000^{2\frac{j}{99}}, j = 0,\dots, 99 \right\} \text{without}\)
\(\text{replacement, and } \alpha_i = 1000 \text{ for } i = 1.\)
\(\text{The local optima } \mathbf{y}_i \text{ are uniformly drawn}\)
\(\text{from the domain } [-5,5]^n \text{ for }\)
\(i = 2,...,101 \text{ and } \mathbf{y}_1 \in [-4,4]^n.\)
\(\text{The global optimum is at } \mathbf{x}^{\text{opt}} = \mathbf{y}_1.\)
\\
\hline
Gallagher's Gaussian
21-hi Peaks Function
&
\(f_{22}(\mathbf{x}) = T_{\text{osz}}\left(10 - \max_{i=1}^{21} w_i \exp\left(- \dfrac{1}{2n} (\mathbf{z} - \mathbf{y}_i)^T \mathbf{B}^T\mathbf{C_i}\mathbf{B} (\mathbf{z} - \mathbf{y}_i) \right) \right)^2\)
\(+ f_{pen}(\mathbf{x}) + \mathbf{f}_{\text{opt}}\)
&
\(w_i = \begin{cases} 1.1 + 8 \times \dfrac{i-2}{19} & \text{for } 2 \leq i \leq 21\\ 10 & \text{for } i = 1 \end{cases}\)
\(\mathbf{B} \text{ is a block-diagonal matrix without}\)
\(\text{permuations of the variables.}\)
\(\mathbf{C_i} = \Lambda^{\alpha_i}/\alpha_i^{1/4} \text{where } \Lambda^{\alpha_i} \text{ is defined as usual,}\)
\(\text{but with randomly permuted diagonal elements.}\)
\(\text{For } i=2,\dots, 21, \alpha_i \text{ is drawn uniformly}\)
\(\text{from the set } \left\{1000^{2\frac{j}{19}}, j = 0,\dots, 19 \right\} \text{without}\)
\(\text{replacement, and } \alpha_i = 1000^2 \text{ for } i = 1.\)
\(\text{The local optima } \mathbf{y}_i \text{ are uniformly drawn}\)
\(\text{from the domain } [-4.9,4.9]^n \text{ for }\)
\(i = 2,...,21 \text{ and } \mathbf{y}_1 \in [-3.92,3.92]^n.\)
\(\text{The global optimum is at } \mathbf{x}^{\text{opt}} = \mathbf{y}_1.\)
\\
\hline
Katsuura Function
&
\(f_{23}(\mathbf{x}) = \left(\dfrac{10}{n^2} \prod_{i=1}^{n} \left( 1 + i \sum_{j=1}^{32} \dfrac{\vert 2^j z_i - [2^j z_i]\vert}{2^j}\right)^{10/n^{1.2}} - \dfrac{10}{n^2}\right)\)
\(+ f_{pen}(\mathbf{x}) + \mathbf{f}_{\text{opt}}\)
&
\(\mathbf{z} = \mathbf{Q}\mathbf{\Lambda}^{100} \mathbf{R}(\mathbf{x} - \mathbf{x}^{\text{opt}})\)
\(\text{ with } \mathbf{R} = P_{11}B_1P_{12}, \mathbf{Q} = P_{21}B_2P_{22}\)
\\
\hline
Lunacek bi-Rastrigin Function
&
\(f_{24}(\mathbf{x}) = \gamma(n)\times\Big(\min\big( \sum_{i=1}^{n} (\hat{x}_i - \mu_0)^2, n + s\sum_{i=1}^{n}(\hat{x}_i - \mu_1)^2\big)\)
\(+ 10 \big(n - \sum_{i=1}^{n}\cos(2\pi z_i) \big)\Big) + 10^{4}f_{pen}(\mathbf{x}) + \mathbf{f}_{\text{opt}}\)
&
\(\mathbf{\hat{x}} = 2 \text{sign}(\mathbf{x}^{\text{opt}}) \otimes \mathbf{x}, \mathbf{x}^{\text{opt}} = 0.5 \mu_0 \mathbf{1}_{-}^{+},\)
\(\mathbf{z} = \mathbf{Q}\mathbf{\Lambda}^{100}\mathbf{R}(\mathbf{\hat{x}} - \mu_0\mathbf{1})\)
\(\text{ with } \mathbf{R} = P_{11}B_1P_{12}, \mathbf{Q} = P_{21}B_2P_{22},\)
\(\mu_0 = 2.5, \mu_1 = -\sqrt{\dfrac{\mu_0^{2} - 1}{s}},\)
\(s = 1 - \dfrac{1}{2\sqrt{n + 20} - 8.2}\)
\\
\hline\end{tabulary}

\end{sidewaystable}\section*{Acknowledgments}
This work was supported by the grant ANR-12-MONU-0009 (NumBBO)
of the French National Research Agency.
This work was further supported by a public grant as part of the Investissement d'avenir project, reference ANR-11-LABX-0056-LMH, LabEx LMH, in a joint call with Gaspard Monge Program for optimization, operations research and their interactions with data sciences.

\begin{sphinxthebibliography}{HAN2016perf}
\bibitem[AIT2016]{\detokenize{AIT2016}}{\phantomsection\label{\detokenize{index:ait2016}} 
O. Ait Elhara, A. Auger, N. Hansen (2016). \href{https://hal.inria.fr/hal-01308566}{Permuted Orthogonal Block-Diagonal
Transformation Matrices for Large Scale Optimization Benchmarking}. GECCO 2016, Jul 2016, Denver,
United States.
}
\bibitem[HAN2009]{\detokenize{HAN2009}}{\phantomsection\label{\detokenize{index:han2009}} 
N. Hansen, S. Finck, R. Ros, and A. Auger (2009).
\href{http://coco.gforge.inria.fr/}{Real-parameter black-box optimization benchmarking 2009: Noiseless
functions definitions}. \href{https://hal.inria.fr/inria-00362633}{Research Report RR-6829}, Inria, updated
February 2010.
}
\bibitem[HAN2016ex]{\detokenize{HAN2016ex}}{\phantomsection\label{\detokenize{index:han2016ex}} 
N. Hansen, T. Tusar, A. Auger, D. Brockhoff, O. Mersmann (2016).
\href{http://numbbo.github.io/coco-doc/experimental-setup/}{COCO: The Experimental Procedure}, \sphinxstyleemphasis{ArXiv e-prints}, \href{http://arxiv.org/abs/1603.08776}{arXiv:1603.08776}.
}
\bibitem[HAN2016perf]{\detokenize{HAN2016perf}}{\phantomsection\label{\detokenize{index:han2016perf}} 
N. Hansen, A. Auger, D. Brockhoff, D. Tusar, T. Tusar (2016).
\href{http://numbbo.github.io/coco-doc/perf-assessment}{COCO: Performance Assessment}. \sphinxstyleemphasis{ArXiv e-prints}, \href{http://arxiv.org/abs/1605.03560}{arXiv:1605.03560}.
}
\bibitem[HAN2016co]{\detokenize{HAN2016co}}{\phantomsection\label{\detokenize{index:han2016co}} 
Nikolaus Hansen, Anne Auger, Olaf Mersmann, Tea Tušar, and Dimo Brockhoff (2016).
\href{http://numbbo.github.io/coco-doc/}{COCO: A Platform for Comparing Continuous Optimizers in a Black-Box
Setting}, \sphinxstyleemphasis{ArXiv e-prints}, \href{http://arxiv.org/abs/1603.08785}{arXiv:1603.08785}.
}
\end{sphinxthebibliography}

\renewcommand{\indexname}{Index}
\printindex
\end{document}